\documentclass[11pt]{article}

\usepackage{amsfonts}
\usepackage{amsmath}

\newcommand{\beq}{\begin{equation}}
\newcommand{\eeq}{\end{equation}}
\newcommand{\bea}{\begin{eqnarray}}
\newcommand{\eea}{\end{eqnarray}}
\newcommand{\beas}{\begin{eqnarray*}}
\newcommand{\eeas}{\end{eqnarray*}}

\newcommand{\acc}{^}

\newcommand{\e}{\varepsilon}

\newtheorem{theorem}{Theorem}[section]

\newtheorem{proposition}[theorem]{Proposition}
\newtheorem{corollary}[theorem]{Corollary}
\newtheorem{lemma}[theorem]{Lemma}
\newtheorem{remark}[theorem]{Remark}
\newtheorem{example}[theorem]{Example}
\newtheorem{examples}[theorem]{Examples}
\newtheorem{foo}[theorem]{Remarks}

\newenvironment{proof}{\addvspace{\medskipamount}\par\noindent{\it
Proof}.}
{\unskip\nobreak\hfill$\Box$\par\addvspace{\medskipamount}}

\title{The subelliptic heat kernel on $\mathbf{SU}(2)$: Representations, Asymptotics and Gradient bounds}

\author{Fabrice Baudoin \footnote{fbaudoin@math.ups-tlse.fr}, Michel Bonnefont \footnote{bonnefon@math.ups-tlse.fr} \\
{\small Institut de Math\'ematiques de Toulouse} \\
{\small Universit\'e de Toulouse} \\
{\small CNRS 5219} \\
}

\begin{document}

\maketitle

\begin{abstract}

The Lie group $\mathbf{SU}(2)$ endowed with its canonical subriemannian structure appears as a three-dimensional 
 model of a \textit{positively curved} subelliptic space. The goal of this work is to study 
the subelliptic heat kernel on it and some related functional inequalities. 

\


\end{abstract}

\tableofcontents

\section{Introduction}

The goal of this work is to study in details the heat kernel and some related functional inequalities in one of the simplest sub-elliptic models
 after the Heisenberg group: the Lie group $\mathbf{SU} (2)$ endowed with its canonical subriemannian structure (coming from
the Hopf fibration $\mathbb{S}^2 \rightarrow \mathbb{S}^3$, see \cite{montgomery}). 
In the classification of three-dimensional homogeneous subriemannian structures (see page 22 in \cite{Vershik}) 
the role played by this  group could be compared to the role played by the sphere in Riemannian geometry: It should be 
a three-dimensional model of a compact positively curved subriemannian space.   

\

In the flat three-dimensional subelliptic model, that is the Heisenberg group, 
the subelliptic heat kernel is quite  well understood. In his celebrated paper \cite{Gaveau}, 
Gaveau gave a useful integral representation and deduced from it small times asymptotics. Since,
 numerous papers have been devoted to the study of this kernel (see for instance \cite{Gaveau2},
 \cite{HQLi} and the references therein). In the case of  $\mathbf{SU} (2)$, we will see that a
 quite similar study can be made: we will obtain an integral representation of the heat kernel
 and will deduce from it the small times asymptotics. These asymptotics give, in particular, 
a way to compute explicitely the Carnot-Carath\'eodory distance associated to the subriemannian 
structure of $\mathbf{SU}(2)$. An interesting fact, on which we will focus, 
is that the Heisenberg group is the tangent cone (in Gromov-Hausdorff sense) to $\mathbf{SU}(2)$ 
so that \textit{dilating} $\mathbf{SU}(2)$, will allow us to recover the known results of the Heisenberg group.

\

On the other hand, recent works have started to study  gradient estimates for subelliptic semigroups
 (see for instance  \cite{bakry-baudoin-bonnefont-chafai}, \cite{DriverMelcher}, \cite{HQLi}, \cite{melcher}). From the point of view of partial differential equations (see \cite{bakry-tata}, \cite{ledoux-zurich}), 
gradient estimates had proved to be a very efficient tool 
for the control of the rate of convergence to equilibrium, quantitative estimates on the
regularization properties of heat kernels, functional inequalities such as
Poincar\'e, logarithmic Sobolev, Gaussian isoperimetric inequalities for heat
kernel measures, etc... When dealing with linear heat equations, those gradient estimates often
rely on the control of the intrinsic Ricci curvature associated to the
generator of the heat equation (Bakry-Emery criterion, see \cite{bakry-emery}). Those methods basically require some form of
ellipticity of the generator and fail in typical subelliptic situations, like for instance in the 
Heisenberg group  (see \cite{bakry-baudoin-bonnefont-chafai}). From the point of view of geometry, these gradient estimates are interesting, because they should contain informations on the curvature of the space. For instance, in Riemannian geometry (see \cite{bakry-st-flour}, \cite{sturm} ), the functional inequality $\| \nabla e^{t \Delta} f \|^2 \le e^{-2 \rho t} e^{t \Delta} (\| \nabla f \|^2)$ is equivalent
to the lower bound $\mathbf{Ricc} \ge \rho$, where $\mathbf{Ricc}$ denotes the Ricci curvature.  In subriemannian 
geometry there is no real analogue of Ricci curvature; for instance, in Lott-Villani-St\"urm sense 
(see \cite{villani-lott}, \cite{sturm1}, \cite{sturm2}), the Ricci curvature of the simplest subelliptic model, the Heisenberg group, is $-\infty$ (see \cite{juillet}). However, we will show in this paper that we obtain exponential decays for the long-time behaviour of gradient estimates of the subelliptic semigroup on the model space  $\mathbf{SU} (2)$ and controls on the small-time behaviour. Nevertheless, as it appears from our methods, the exponential decays we obtain are optimal but are mainly consequences of spectral properties, so that we do not really rely on any notion of \textit{intrinsic Ricci curvature} excepted in the Li-Yau type estimate that we obtain. In the future, we hope to extend those methods to cover more general situations and to make the link with more geometrically oriented works like for instance \cite{rumin}, where a Bonnet-Myers type theorem is obtained in a hypoelliptic situation.

\

So, finally, this work is mainly divided into two parts. In a first part (Section 3), we will study the subelliptic heat kernel on  $\mathbf{SU} (2)$.
We provide its spectral decomposition, prove an integral representation of it and compute its small times asymptotics. 
In the second part (Section 4) we will focus on gradient estimates, using the previous results.

\section{Preliminaries on $\mathbf{SU}(2)$}

In what follows, we consider the Lie group $\mathbf{SU} (2)$, i.e. the group of
$2 \times 2$, complex, unitary matrices of determinant $1$. Its Lie algebra $\mathfrak{su} (2)$ consists of $2 \times 2$, complex,
skew-adjoint matrices of trace $0$. A basis of $\mathfrak{su} (2)$
is formed by the Pauli matrices:
 \[
\text{ }X=\left(
\begin{array}{cc}
~0~ & ~1~ \\
-1~& ~0~
\end{array}
\right) ,\text{ }Y=\left(
\begin{array}{cc}
~0~ & ~i~ \\
~i~ & ~0~
\end{array}
\right),
Z=\left(
\begin{array}{cc}
~i~ & ~0~ \\
~0~ & -i~
\end{array}
\right) ,
\]
for which the following relationships hold
\begin{align}\label{Liestructure}
[Z,X]=2Y, \quad [X,Y]=2Z, \quad [Y,Z]=2X.
\end{align}
We denote $X,Y,Z$ the left invariant vector fields on 
$\mathbf{SU} (2)$ corresponding to the Pauli matrices. 
The Laplace-Beltrami operator for the bi-invariant Riemannian structure of 
$\mathbf{SU} (2) \simeq \mathbb{S}^3$ is
\begin{align*}
\Delta  =X^2+Y^2+Z^2.
\end{align*}
It is in the center of the universal enveloping algebra of the vector fields $X,Y,Z$. In the sequel, we shall mainly 
be interested in the operator
\[
\mathcal{L} =X^2+Y^2.
\]

According to the relations (\ref{Liestructure}) and due to H\"ormander's theorem, 
$\mathcal{L}$  is subelliptic but not elliptic so that the associated geometry is not Riemannian 
but only subriemannian.

Associated to $\mathcal{L}$, there is a notion of length of gradient that is given via
the \textit{carr\'e du champ} operator defined for smooth functions by
\[
\Gamma (f,f)=\frac{1}{2} (\mathcal{L}f^2-2f\mathcal{L}f)=(Xf)^2+(Yf)^2.
\] 
The intrinsic distance associated to $\mathcal{L}$  is given

\[
\delta(g_1,g_2)=\sup_{f \in \mathcal{C}} \{ \mid f(g_1) -f(g_2) \mid \}
\]

where $\mathcal{C}$ is the set of smooth maps $\mathbf{SU}(2) \rightarrow \mathbb{R}$ that satisfy
 $\Gamma (f,f) (x) \le 1$ , $ x \in  \mathbf{SU}(2)$. This distance is the 
Carnot-Carath\'eodory distance. Via Chow's theorem, it can also be defined as the minimal length
of horizontal curves joining two given points (see Chapter 3 of \cite{baudoin}).

\

To study $\mathcal{L}$, we will use the cylindric coordinates introduced in \cite{Cow-Sik}:
\[
(r,\theta, z) \rightarrow \exp \left(r \cos \theta X +r \sin \theta Y \right) \exp (z Z)= \left(
\begin{array}{cc}
 \cos(r) e\acc{iz}& \sin(r)e\acc{i(\theta -z)}\\
 -\sin(r) e\acc{-i(\theta -z)}& \cos(r) e\acc{-iz}
\end{array}\right),
\]
with 
\[
0 \le r < \frac{\pi}{2}, \text{ }\theta \in [0,2\pi], \text{ }z\in [-\pi,\pi].
\]

Simple but tedious computations show that in these coordinates, the left-regular representation sends the matrices
 $X$, $Y$ and $Z$ to the
left-invariant vector fields: 
\[
X=\cos (-\theta +2z) \frac{\partial}{\partial r}+\sin (-\theta +2z) \left( \tan r \frac{\partial}{\partial z}
+\left(\tan r+\frac{1}{\tan r}\right)  \frac{\partial}{\partial \theta}\right),
\]
\[
Y=-\sin(2z-\theta) \frac{\partial}{\partial r}+\cos (2z-\theta) \left( \tan r \frac{\partial}{\partial z}
+\left(\tan r+\frac{1}{\tan r}\right)  \frac{\partial}{\partial \theta}\right),
\]
\[
Z=\frac{\partial}{\partial z},
\]
and that the bi-invariant normalized Haar measure reads:
\[
d\mu=\frac{1}{4 \pi^2} \sin 2r dr d\theta dz
\]

\begin{remark}
The right regular representation sends the matrices $X$, $Y$ and $Z$ to the
right-invariant vector fields
\[
\hat{X}=\cos \theta  \frac{\partial}{\partial r}+\sin \theta  \left( \tan r \frac{\partial}{\partial z}
+\left(\tan r-\frac{1}{\tan r}\right)  \frac{\partial}{\partial \theta}\right)
\]
\[
\hat{Y}=\sin \theta \frac{\partial}{\partial r}-\cos \theta \left( \tan r \frac{\partial}{\partial z}
+\left(\tan r-\frac{1}{\tan r}\right)  \frac{\partial}{\partial \theta}\right).
\]
\[
\hat{Z}=\frac{\partial}{\partial z}+2 \frac{\partial}{\partial \theta}
\]
\end{remark}

We therefore obtain
\begin{align*}
\mathcal{L} & =X^2+Y^2 \\
 & =\frac{\partial^2}{\partial r^2}+2 \text{ cotan } 2r\frac{\partial}{\partial r}+\left( 2+\frac{1}{\tan^2 r}
+\tan^2 r  \right) \frac{\partial^2}{\partial \theta^2} +\tan^2 r \frac{\partial^2}{\partial z^2}+2(1+\tan^2r)
\frac{\partial^2}{\partial z \partial \theta}
\end{align*}
and
\begin{align*}
\Delta & =X^2+Y^2+Z^2 \\
 & =\frac{\partial^2}{\partial z^2} +\mathcal{L} 
\end{align*}

Note that $\mathcal{L}$ commutes with  $\frac{\partial}{\partial \theta}$ and with  $\frac{\partial}{\partial z}$.

\begin{remark}[Probabilistic interpretation]
The computation of the left-regular representation shows that if $(X_t)_{t \ge 0}$ is the Markov process that is the matrix-valued solution of the stochastic differential equation (written in Stratonovitch form)
\[
dX_t =X_t (X \circ dB^1_t + Y \circ dB^2_t), \quad X_0=1,
\]
where $(B^1_t,B^2_t)_{t \ge 0}$ is a two-dimensional Brownian motion, then, in law,
\[
X_t=\exp \left( \rho_t (X\cos \theta_t  + Y\sin \theta_t )\right) \exp (z_t Z), \quad t\ge 0,
\]
where $(\rho_t,\theta_t, z_t)_{t \ge 0}$ solve the following stochastic differential equations (written in It\^o's form):
\[
d\rho_t =2 \mathrm{ cotan} 2 \rho_t dt +dB^1_t,
\]
\[
d\theta_t=\frac{2}{\sin \theta_t} dB^2_t,
\]
\[
dz_t =\tan \rho_t dB^2_t.
\]
\end{remark}

\section{The subelliptic heat kernel on $\mathbf{SU}(2)$}

By hypoellipticity, the heat semigroup $P_t=e^{t \mathcal{L}}$ admits
a smooth kernel with respect to the Haar measure $\mu$ of $\mathbf{SU} (2)$. Our goal in this section will be to 
derive various representations of this kernel and to get precise asymptotics in small times.

\subsection{Spectral decomposition of the heat kernel}

Since $\mathcal{L}$ commutes with $\frac{\partial}{\partial \theta}$ that vanishes at 0, we deduce that the heat
 kernel (issued from the identity)
 of $P_t=e^{t\mathcal{L}}$, only depends on $(r,z)$. It will be denoted by $p_t (r,z)$.

We first obtain the spectral decomposition of $p_t (r,z)$:

\begin{proposition}\label{spectral_decomposition}
For $t >0$, $0 \le r < \frac{\pi}{2}, \text{ }z\in [-\pi,\pi]$.
\[
p_t (r,z)=\sum_{n=-\infty}^{+\infty}\sum_{k=0}^{+\infty}(2k+\mid n \mid +1)e^{-(4k(k+\mid n \mid +1)+2\mid n \mid  )t}
 e^{inz} (\cos r )^{\mid n \mid}  P_k^{0,\mid n \mid} (\cos 2r ),
\]
where
$$P_k^{0,\mid n \mid} (x)=\frac{(-1)^k}{2^k k ! (1+x)^{\mid n \mid} }\frac{d^k}{dx^k} \left( (1+x)^{\mid n \mid} (1-x^2)^k \right)
$$
is a Jacobi polynomial.
\end{proposition}

\begin{proof}
Since the points $(r,\theta, z)$ and $(r,\theta , z+2 \pi)$ are the same, we can define $p_t(r,z)$ for all $z\in \mathbb R$ and it is $2\pi$-periodic.
The idea is then to expand $p_t (r,z)$ as a Fourier series in $z$:
\[
p_t (r,z)=\sum_{n=-\infty}^{+\infty}e^{inz} \Phi_n (t,r)
\]
Since $p_t(r,z)$ satisfies the partial differential equation,
$$\frac{\partial p_t}{\partial t}=\mathcal L p_t,
$$
we obtain for $\Phi_n$ the following equation
\[
\frac{\partial \Phi_n}{\partial t}=\frac{\partial^2 \Phi_n}{\partial r^2}
+2 \text{ cotan } 2r\frac{\partial \Phi_n}{\partial r}-n^2 \tan^2 r \Phi_n 
\]
and look for a solution under the form
\[
\Phi_n (t,r)=e^{-2nt}(\cos r)^{\mid n \mid} g_n (t,\cos 2r).
\]
We get: 
\[
\frac{\partial g_n}{\partial t}=4\mathcal G _n (g_n)
\] 
where 
\[
\mathcal G_n=
(1-x\acc2)\frac{\partial^2 }{\partial x^2}+ ( |n|-(2+|n|)x)\frac{\partial }{\partial x} 
\]


It is well-known that eigenvectors of $\mathcal G_n$ are the Jacobi polynomials:
$$P_k^{0,\mid n \mid} (x)=\frac{(-1)^k}{2^k k ! (1+x)^{\mid n \mid} }\frac{d^k}{dx^k} \left( (1+x)^{\mid n \mid} (1-x^2)^k \right)
$$
In fact we have 
$$\mathcal G_n (P_k\acc{0 ,|n|}) (x)= -k(k+n+1) P_k\acc{0,|n|}(x)
$$
So we are finally led to put
$$p_t (r,z)=\sum_{n=-\infty}^{+\infty}\sum_{k=0}^{+\infty} \alpha_{k,n}e^{-(4k(k+\mid n \mid +1)+2\mid n \mid  )t} e^{inz} (\cos r )^{\mid n \mid}  P_k^{0,\mid n \mid} (\cos 2r )
$$ for some $\alpha_{k,n}$, where the $\alpha_{k,n}$ will be determined by the initial condition at time 0.

Clearly $p_t$ satisfies the equation $\frac{\partial p_t}{\partial t}=\mathcal L p_t$ and, by using the fact that
$(P_k\acc{0,|n|})_{k\geq 0}$ is an orthogonal
 basis of ${L\acc 2 ([-1,1],(1+x)\acc {|n|} du)}$ with ${||P_k\acc{0,|n|}||\acc2= \frac{2\acc{|n|+1}}{2k+|n|+1}}$ 
we easily check that for a smooth $f$
\[
\frac{1}{2\pi}\int_0\acc \frac{\pi}{2}\int_0\acc {2\pi} p_t(r,z) f(r,z)\sin(2r)dr dz 
\rightarrow_{t\rightarrow 0} f(0,0) 
\]
as soon as $\alpha_{k,n}=2k+|n|+1$.

\end{proof}

\begin{remark}
By using the representation theory of $\mathbf{SU}(2)$, a similar spectral decomposition is given in \cite{bauer}. 
Nevertheless, for the sake of completeness, we included this elementary proof.
\end{remark}

\subsection{Integral representation of the heat kernel}

We now provide an integral representation of $p_t$  based on the following formula:
\[
e^{t\mathcal{L}}=e^{-t\frac{\partial^2}{\partial z^2}} e^{t\Delta}
\]
that stems from the commutation between $\Delta$ and $\frac{\partial}{\partial z}$. Since $\Delta$ is the
Laplace-Beltrami operator on the three-dimensional sphere which has a well-known  heat kernel, it will lead  
to an expression of $p_t$.

Let us consider on the interval $[-1,1]$ the second order differential operator
\[
\mathcal J =(1-x^2)\frac{d^2}{dx^2}-3x \frac{d}{dx}.
\]
For $m \ge 0$, let $U_m$ denotes the Chebyshev polynomial of the second kind:
\[
U_m (\cos x)=\frac{\sin(m+1)x}{\sin x},
\]
and
\begin{align}\label{spectralQ}
q_t (x)=\sum_{m=0}^{+\infty} (m+1) e^{-m(m+2)t} U_m (x), \quad x  \in [-1,1].
\end{align}

It is known that if $f$ is a smooth function $[-1,1] \rightarrow \mathbb{R}$, then
\[
(e^{t \mathcal{J}}f)(1)=\frac{2}{\pi} \int_{-1}^1 q_t (x) f(x)(1-x^2)^{1/2} dx,
\]
\begin{lemma}
If $f$ is a smooth function $\mathbf{SU}(2)  \rightarrow \mathbb{R}$, then for $t \ge 0$,
\[
(e^{t \Delta}f)(0)
=\frac{1}{4\pi^2}\int_0^{\frac{\pi}{2}} \int_0^{2\pi} \int_{-\pi}^{\pi} 
q_t (\cos r \cos z) f(r,\theta,z) \sin 2r dr d\theta dz
\]
\end{lemma}
\begin{proof} 
An easy calculation shows that the function $q_t(\cos r \cos z)$ solves the heat equation
\[
\frac{\partial }{\partial t} (q_t(\cos r \cos z))  = \Delta (q_t(\cos r \cos z)).
\]
Now we have to check the initial condition. We must show 
$$\frac{1}{2\pi}\int_{0}\acc{\pi/2}\int_{0}^{2\pi}q_t(\cos r \cos z)f(r,z)\sin 2r dr dz \longrightarrow _{t\rightarrow 0} f(0,0)
$$
Since we will make the following change of variables:
$$\left\{\begin{array}{ccc}
u&=& \cos r\cos z\\
v&=&\cos r \sin z
\end{array}\right.
$$
we take the function $f$ of the form $f(r,z)=g(\cos r \cos z)h(\cos r \sin z)$. The new domain is $D=\{(x,y), x\acc 2+y\acc 2 \leq 1\}$ and the Jacobian determinant is $\frac{1}{2}\sin 2r$. So
\begin{eqnarray*}
& &\frac{1}{2\pi}\int_{r=0}\acc{\pi/2}\int_{z=0}\acc{2\pi}q_t(\cos r \cos z)g(\cos r \cos z)h(\cos r \sin z)\sin 2r dr dz\\
&=& \frac{1}{\pi}\int \int_D q_t(u)g(u)h(v)du dv\\
&= &\frac{1}{\pi}\int_{-1}\acc 1 \left(\int_{-(1-u\acc 2)\acc {1/2}}\acc{(1-u\acc 2)\acc {1/2}} h(v) dv\right)q_t(u)g(u) du
\end{eqnarray*}

We may rewrite it as
$$\frac{2}{\pi}\int_{-1}\acc 1 q_t(u)l(u)(1-u\acc 2)\acc {1/2}du 
$$
where $ l$ is the continuous fonction 
$$l(u)=g(u) \left( \frac{\int_{-(1-u\acc 2)\acc {1/2}}\acc{(1-u\acc 2)\acc {1/2}} h(v) dv}{2(1-u\acc 2)\acc {1/2}}\right)
$$
Now, since $q_t$ is the heat kernel of a diffusion issued of $1$ with respect to the measure 
$\frac{2}{\pi} (1-u\acc 2)\acc{1/2} du$ and $l$ is continuous, the last quantity is converging towards $l(1)=g(1)h(0)=f(0,0)$ and the lemma is proved.
\end{proof}

\begin{remark}
The previous lemma shows that if $\rho$ is the Riemannian distance from 0, then in our cylindric coordinates, we have
\[
\cos \rho =\cos r \cos z.
\]
\end{remark}

From the previous proposition, we can now derive  an expression for $p_t$ in terms of $q_t$.

Let us first describe some properties of $q_t$ that will be useful in the sequel.
From the Poisson summation formula, we obtain that for $\theta \in \mathbb{R}$:
\begin{align*}
q_t (\cos \theta)&=\frac{\sqrt{\pi}e^t}{4 t^{\frac{3}{2}}}  \frac{1}{\sin \theta} \sum_{k \in \mathbb{Z}} (\theta +2k\pi)e^{-\frac{(\theta+2k\pi)^2}{4t}}\\
 & =\frac{\sqrt{\pi}e^t}{4 t^{\frac{3}{2}}} \frac{\theta}{\sin \theta} e^{-\frac{\theta^2}{4t}} \left(1 +2\sum_{k=1}^{+\infty} e^{-\frac{k^2\pi^2}{t}}\left( \cosh \frac{k\pi \theta}{t} +2k\pi \frac{\sinh \frac{k\pi \theta}{t}}{\theta} \right) \right)
\end{align*}
These expressions show that $q_t (\cos \theta)$ admits an analytic extension for $\theta \in \mathbb{C}$. We moreover obtain precise estimates:

\
\begin{itemize}
\item Let $\varepsilon>0$, for $x \in (-1+\varepsilon,1]$ and $t >0$:
\begin{align}\label{estimate1}
q_t (x)=\frac{\sqrt{\pi}e^t}{4 t^{\frac{3}{2}}}  \frac{\mathrm{arcos }x}{\sqrt{1-x^2}} e^{-\frac{(\mathrm{arcos }x)^2}{4t}} \left(1 +R_1(t,x) \right),
\end{align}
where for some positive constants $C_1$ and $C_2$ depending only in $\varepsilon$, $\mid R_1(t,x) \mid \le C_1 e^{-\frac{C_2}{t}}$.
\item For $x \in [1,+\infty)$ and $t >0$:
\begin{align}\label{estimate2}
q_t (x)=\frac{\sqrt{\pi}e^t}{4 t^{\frac{3}{2}}} \frac{\mathrm{arcosh }x}{\sqrt{x^2-1}} e^{\frac{(\mathrm{arcosh }x)^2}{4t}} \left(1 +R_2(t,x) \right),
\end{align}
where for some positive constants $C_3$ and $C_4$, $\mid R_2(t,x) \mid \le C_3 e^{-\frac{C_4}{t}}$.
\end{itemize}

\begin{proposition}\label{representation}
We have for $t >0$, $r \in [0,\pi /2)$, $z \in [-\pi, \pi]$,
\[
p_t (r,z)=\frac{1}{\sqrt{4\pi t}} \int_{-\infty}^{+\infty} e^{-\frac{(y+iz)^2}{4t}} q_t (\cos r \cosh y) dy
\]
\end{proposition} 
\begin{proof}
Let
\[
h_t (r,z)=\frac{1}{\sqrt{4\pi t}} \int_{-\infty}^{+\infty} e^{-\frac{(y+iz)^2}{4t}} q_t (\cos r \cosh y) dy;
\]
the integral being well defined thanks to the estimates on $q_t$.
By using the fact that
\[
\frac{\partial}{\partial t} \left(\frac{ e^{-\frac{y^2}{4t}}}{\sqrt{4\pi t}}\right)=\frac{\partial^2}{\partial y^2} \left(\frac{ e^{-\frac{y^2}{4t}}}{\sqrt{4\pi t}}\right)
\]
and
\[
\frac{\partial }{\partial t} (q_t(\cos r \cos z))  = \Delta (q_t(\cos r \cos z)),
\]
a double integration by parts with respect to the variable $y$ shows that 
\[
\frac{\partial h_t}{\partial t}=\mathcal L h_t.
\]
Let us now check the initial condition. 
Let $f(r,z)=e^{i\lambda z} g(r)$ where $\lambda \in \mathbb{R}$ and $g$ is a smooth function.
We have
\[
\frac{1}{4\pi^2}\int_0^{\frac{\pi}{2}} \int_0^{2\pi} \int_{-\pi}^{\pi} 
h_t (r , z) f(r,z) \sin 2r dr d\theta dz=e^{t\lambda^2} \left( e^{t \Delta} g \right) (0),
\]
so that we obtain the required result.
\end{proof}
We are now in position to collect some properties of $p_t$.
\begin{proposition}\label{Laplacetransform}
For $\lambda \in \mathbb{C}$, $\mathbf{Re} \lambda >0$, $r \in [0,\pi /2)$, $z \in [-\pi, \pi]$,

\[
\int_0^{+\infty}  p_{t} (r,z) e^{-t-\frac{\lambda}{t}} dt=\int_{-\infty}^{+\infty}\frac{dy}{8 \pi^2\left( \cosh \sqrt{y^2+4\lambda} -\cos r \cos (z+iy) \right) }
\]

\end{proposition}

\begin{proof}
We have
\begin{align*}
 \int_0^{+\infty} e^{-\frac{\lambda}{t}} p_{t} (r,z) e^{-t} dt = \frac{1}{\sqrt{4\pi}}\int_{-\infty}^{+\infty} \int_{0}^{+\infty} e^{-t -\frac{y^2+4\lambda}{4t}}
 q_t (\cos r \cos (z+iy) ) \frac{dt}{\sqrt{t}} dy
\end{align*}
We now compute
\[
\int_0^{+\infty} e^{-t -\frac{y^2+4\lambda}{4t}}
 q_t (\cos r \cos (z+iy) ) \frac{dt}{\sqrt{t}} 
\]
by using the symbolic calculus on differential operators (it can be made rigorous with \ref{spectralQ}).
\begin{align*}
\int_0^{+\infty} e^{-t -\frac{y^2+4\lambda}{4t}}
 e^{t\Delta} \frac{dt}{\sqrt{t}} dt&=\int_0^{+\infty} e^{-\frac{y^2+4\lambda}{4t}}
 e^{-t(-\Delta+1)} \frac{dt}{\sqrt{t}} dt \\
 &= \frac{\sqrt{\pi}}{\sqrt{-\Delta +1}}e^{-\sqrt{y^2+4\lambda} \sqrt{-\Delta+1}}
\end{align*}
But from Taylor \cite{taylor2} pp. 95,
\[
 \frac{1}{\sqrt{-\Delta +1}}e^{-\sqrt{y^2+4\lambda} \sqrt{-\Delta+1}}=
\frac{1}{4 \pi^2\left( \cosh \sqrt{y^2+4\lambda} -\cos r \cos z \right) },
\]
which implies the result.
\end{proof}
If we fix, $r \in [0,\pi /2)$, $z \in [-\pi, \pi]$, we observe  that it possible to find $\theta (r,z) \in \mathbb{R}$, such that for $\lambda \in \mathbb{C},  \mathbf{Re} \lambda >0 $ and $y \in \mathbb{R}$,
\[
\cosh \sqrt{y^2+4\lambda} = \cos r \cos (z+iy) \Rightarrow \mathbf{Re} \lambda \le \theta (r,z),
\]
where we use the principal branch of the square root.   
By inverting the last Laplace transform of the previous proposition, we therefore get:
\begin{corollary}
We have for $t >0$, $r \in [0,\pi /2)$, $z \in [-\pi, \pi]$, and $\gamma > \theta (r,z)$,
\[
p_t (r,z)=\frac{e^t}{16 i\pi^3 t^2} \int_{\lambda=\gamma-i\infty}^{\gamma+i\infty}\int_{y=-\infty}^{+\infty}  
 \frac{e^{\frac{\lambda}{t}}dy d\lambda }{\cosh \sqrt{y^2+4\lambda}-\cos r \cos(z+iy)}
\]
\end{corollary}

>From Proposition \ref{Laplacetransform}, we also deduce:
\begin{proposition}
The Green function of the operator $-\mathcal{L} +1$ is given by
\[
G(r,z)=\frac{1}{8 \pi} \frac{1}{\sqrt{1-2\cos r \cos z +\cos^2 r } }.
\]
\end{proposition}

\begin{proof}
Let us assume $r \neq 0$, $z \neq 0$. In that case the Laplace transform of Proposition \ref{Laplacetransform} can be 
extended to $\lambda =0$ and we have:
\begin{align*}
G(r,z)& =\int_{-\infty}^{+\infty}\frac{dy}{8 \pi^2\left( \cosh y -\cos r \cos (z+iy) \right) } \\
  &=\int_{-\infty}^{+\infty}\frac{dy}{8 \pi^2\left( (1-\cos r \cos z)\cosh y -i \cos r \sin z \sinh y \right) } \\
  &=\frac{1}{8 \pi^2} \frac{1}{\sqrt{1-2\cos r \cos z +\cos^2 r }} \int_{-\infty}^{+\infty}\frac{dy}{\cosh y} \\
   & = \frac{1}{8 \pi} \frac{1}{\sqrt{1-2\cos r \cos z +\cos^2 r } }
\end{align*}
\end{proof}
\subsection{Asymptotics of the heat kernel in small times}

The goal of this section is to obtain the precise asymptotics of the heat kernel when $t \to 0$. We start with the points of the form $(0,z)$ that lie on the cut-locus of 0. 

\begin{proposition}\label{asymptotics1}
For $t>0$ and $z \in [0,\pi)$,
\[
p_t (0,z)=\frac{\pi^2 e^t}{4t^2} e^{-\frac{2\pi z -z^2}{4t}}\sum_{k \in 
\mathbb{Z}}e^{-\frac{k(k+1)\pi^2}{t}}
\frac{ (2k+1)+2k e^{-\frac{\pi}{2t} (z+2k\pi)}} { \left( 1+e^{-\frac{\pi}{2t} 
(z+2k\pi)} \right)^2}
\]
therefore, when $t \rightarrow 0$,
\[
p_t(0,z)=\frac{\pi^2 e^t}{4t^2} e^{-\frac{2\pi z -z^2}{4t}} \left( 
1+O(e^{-\frac{C}{t}}) \right)
\]
\end{proposition}

\begin{proof}
Let $z \in (0,\pi]$. We have
\[
p_t (0,z)=\frac{1}{\sqrt{4\pi t}} \int_{-\infty}^{+\infty} e^{-\frac{y^2}{4t}} 
q_t ( \cosh (y-iz))) dy,
\]
But
\[
q_t (\cosh (y-iz))=\frac{\sqrt{\pi}e^t}{4 t^{\frac{3}{2}}}  \frac{1}{\sinh (y-iz)}
 \sum_{k \in \mathbb{Z}} (y-iz -2ik\pi)e^{\frac{(y-iz-2ik\pi)^2}{4t}}\\
\]
and for $k\in \mathbb{Z}$, from the residue theorem,
\[
\int_{-\infty}^{+\infty} \frac{y-iz-2ik\pi}{\sinh (y-iz)} e^{-\frac{iy}{2t}(z+2k\pi)}dy
=2\pi^2 e^{\frac{(z+2k\pi)^2-(2k+1)\pi(z+2k\pi) }{2t} }
\frac{ (2k+1)+2k e^{-\frac{\pi}{2t} (z+2k\pi)}} { \left( 1+e^{-\frac{\pi}{2t} 
(z+2k\pi)} \right)^2}.
\]
The result easily follows.
\end{proof}

We now come to points $(r,z)$ that do not lie on the cut-locus, that is $r \neq 0$.

\begin{proposition}
For $r \in \left(0,\frac{\pi}{2}\right)$, when $ t \to 0$,
\[
p_t (r,0) \sim   \frac{r}{\sin r} \sqrt{\frac{1}{1-r \mathrm{cotan} r} }\frac{\sqrt{\pi}e^{-\frac{r^2}{4t}}}{4t^{\frac{3}{2}}}.
\]
\end{proposition}

\begin{proof}
We fix $r \in \left(0,\frac{\pi}{2}\right)$. From the proposition \ref{representation} and due to the estimates on $q_t$ we get:
\begin{align*}
p_t (r,0) \sim_{t \to 0} \frac{1}{8t^2}(J_1 (t)+J_2(t)),
\end{align*}
where
\[
J_1 (t)=\int_{ \cosh y \le \frac{1}{\cos r}} e^{- \frac{y^2 + (\mathrm{arcos} (\cos r \cosh y))^2}{4t}}
\frac{\mathrm{arcos} (\cos r \cosh y)}{\sqrt{1 -\cos^2 r \cosh^2 y}} dy
\]
and
\[
J_2 (t)=\int_{ \cosh y \ge \frac{1}{\cos r}} e^{- \frac{y^2 - (\mathrm{arcosh} (\cos r \cosh y))^2}{4t} }\frac{\mathrm{arcosh} (\cos r \cosh y)}{\sqrt{\cos^2 r \cosh^2 y-1}} dy.
\]
We now analyze the two above integrals in small times thanks to the Laplace method and show that  $J_2(t)$ can be omitted.

On the interval $\left[ -\mathrm{arcosh} \frac{1}{\cos r} ,\mathrm{arcosh} \frac{1}{\cos r} \right]$, the function
\[
f(y)=y^2 + (\mathrm{arcos} (\cos r \cosh y))^2
\]
has a unique minimum which is attained at $y=0$ and, at this point:
\[
f''(0)=2(1-r\mathrm{cotan} r).
\]
Therefore, thanks to the Laplace method
\[
J_1 (t)\sim_{t \to 0} e^{-\frac{r^2}{4t}} \frac{r}{\sin r} \sqrt{\frac{4\pi t}{1-r \mathrm{cotan} r} }.
\]
We now analyze the second integral. On $\left( -\infty , -\mathrm{arcosh} \frac{1}{\cos r} \right) \cup \left( \mathrm{arcosh} \frac{1}{\cos r}, +\infty \right)$, the function
\[
g(y)=y^2 - (\mathrm{arcosh} (\cos r \cosh y))^2,
\]
has no minimum. Therefore, from the Laplace method $J_2 (t)$ is negligeable with respect to $J_1(t)$ when $t \to 0$.
\end{proof}

The previous proposition can be extended by the same method when $z \neq 0$. 
If we fix $r\in \left( 0, \frac{\pi}{2} \right), z \in [-\pi, \pi]$, then the function
\[
f(y)=(y-iz)^2 + (\mathrm{arcos} (\cos r \cosh y))^2,
\]
defined on the strip $\mid \mathrm{Re} (y) \mid  < \mathrm{arcosh} \frac{1}{\cos r}$ has a critical point at $i\theta (r,z)$ where $
\theta(r,z)$ is the unique solution in $[-\pi,\pi]$ to the equation:
\[
\theta (r,z)-z=\cos r \sin \theta (r,z)\frac{ \mathrm{arcos} (\cos \theta (r,z) \cos r ) }{\sqrt{1-\cos^2 r \cos^2 \theta (r,z)}}.
\]
Indeed, with $u= \cos r \cos \theta$
$$\frac{\partial}{\partial \theta} \left( \theta - \cos r \sin \theta \frac{ \mathrm{arcos} (\cos \theta  \cos r ) }{\sqrt{1-\cos^2 r \cos^2 \theta}}\right)=  \frac{\sin^2 r}{1-u(r,z)^2} \left( 1-  \frac{u(r,z) \mathrm{arcos } u(r,z)}{\sqrt{1-u^2(r,z)}} \right)
$$
which is positive. So this last function is bijective from $[-\pi,\pi]$ on itself.

We observe that at the point $\theta (r,z)$,  $f''(i\theta (r,z))$ is a positive real number:

$$ f''(i\theta (r,z))= 2  \frac{\sin^2 r}{1-u(r,z)^2} \left( 1-  \frac{u(r,z) \mathrm{arcos } u(r,z)}{\sqrt{1-u^2(r,z)}} \right)
$$
where $u(r,z)= \cos r \cos \theta (r,z)$.
By the same method than in the previous proposition, we obtain:
\begin{proposition}
Let $r\in \left( 0, \frac{\pi}{2} \right), z \in [-\pi, \pi]$. When $t \to 0$,

$$p_t(r,z) \sim \frac{1}{\sin r} \frac{\mathrm{arccos } u(r,z)}{\sqrt{ 1-  \frac{u(r,z) \mathrm{arcos } u(r,z)}{\sqrt{1-u^2(r,z)} }}} \frac{\sqrt{\pi}e^{-\frac{(\theta (r,z)-z)^2 \tan^2 r}{4t \sin^2 \theta (r,z)}}}{4 t^{\frac{3}{2}}}.
$$
\end{proposition}

\begin{remark}
According to L\'eandre results \cite{leandre1} and \cite{leandre2}, the previous asymptotics give a way to compute the sub-Riemannian distance from 0 to the point $(r,\theta,z) \in \mathbf{SU}(2)$ by computing $\lim_{t \to 0} -4t \ln p_t (r,z)$. This distance does not depend on the variable $\theta$ and shall be denoted by $d(r,z)$.
\begin{itemize}
\item For $z \in [-\pi,\pi]$, 
\[
d^2 (0,z)=2\pi \mid z \mid -z^2 .
\]
\item For $z \in [-\pi,\pi]$, $r\in \left( 0, \frac{\pi}{2} \right)$,  
\[
d^2(r,z)=\frac{(\theta(r,z)-z)^2 \tan^2 r}{ \sin^2 \theta (r,z)}.
\]
In particular, $d^2(r,0)=r^2$.

\end{itemize}
In particular, the sub-Riemannian diameter of $\mathbf{SU}(2)$ is thus $\pi^2$.
\end{remark}

\subsection{From $\mathbf{SU}(2)$ to the Heisenberg group}

Our goal in this section is to exhibit the close connection that exists between the subelliptic operator
$\mathcal{L}$ on $\mathbf{SU}(2)$ and the canonical subelliptic operator on the Heisenberg group $\mathbb{H}$.

Let us first recall some basic properties of the three-dimensional Heisenberg group 
(see by e.g. \cite{baudoin}, \cite{bakry-baudoin-bonnefont-chafai} and the references therein): $\mathbb{H}$ can
be represented as $\mathbb{R}^3$ endowed with the polynomial group law:
\[
(x_1,y_1,z_1) (x_2,y_2,z_2)=(x_1+x_2,y_1+y_2,z_1+z_2+x_1y_2-x_2y_1).
\]
The left invariant vector fields read in cylindric coordinates ($x=r \cos \theta, y=r\sin \theta$):
\[
\tilde{X}= \cos \theta \frac{\partial }{\partial r} -\frac{\sin
\theta}{r}\frac{\partial }{\partial \theta}- r \sin
\theta \frac{\partial }{\partial z}
\]
\[
\tilde{Y}= \sin \theta \frac{\partial }{\partial r} +\frac{\cos
\theta}{r} \frac{\partial }{\partial \theta} + r \cos
\theta \frac{\partial }{\partial z}
\]
\[
\tilde{Z}=\frac{\partial }{\partial z}.
\]
And the following equalities hold
\[
\lbrack \tilde{X},\tilde{Y}]=2\tilde{Z},\text{ }[\tilde{X},\tilde{Z}]=[\tilde{Y},\tilde{Z}]=0.
\]
We denote
\[
\tilde{\mathcal{L}}=\tilde{X}^2 + \tilde{Y}^2.
\]

Due to Gaveau's formula (see \cite{Gaveau}), with respect to the Lebesgue measure $r dr d\theta dz$ the heat kernel associated 
to the semigroup $(\tilde{P}_t)_{t\geq0}=(e^{t\tilde{\mathcal{L}}})_{t\ge0}$ writes
\begin{equation}\label{gaveau}
h_t(r,z)= \frac{1}{16 \pi^2} \int_{-\infty}^{+\infty} e^{\frac{i \lambda z}{2}} \frac{\lambda}{\sinh \lambda t}
e^{-\frac{r^2}{4}\lambda  \text{cotanh}  \lambda  t }d \lambda.
\end{equation}

If $d$ denotes the Carnot Carth\'eodory distance on $\mathbf{SU}(2)$, then it is known that
the Heisenberg group is the tangent cone in Gromov-Hausdorff sense. More precisely from 
Mitchell theorem \cite{Mit} (see also \cite{baudoin}), for any positive $R$:
\[
\lim_{n \rightarrow +\infty} \delta_{\mathbf{GH}} \left(
n\mathbf{B}_{\mathbf{SU}(2)} \left( 0,R \right) ,
\mathbf{B}_{\mathbb{H}} \left( 0 ,R \right) \right)=0,
\]
where:
\begin{itemize}
\item $n\mathbf{B}_{\mathbf{SU}(2)} \left( 0,R \right) $ is the open ball in $\mathbf{SU}(2)$
with radius $R$ for the dilated Carnot-Carath\'eodory metric $nd$;
\item $\mathbf{B}_{\mathbb{H}} \left( 0 ,R \right)$ is the open ball in $\mathbb{H}$
with radius $R$ for the Carnot-Carath\'eodory metric;
\item $ \delta_{\mathbf{GH}}$ is the Gromov-Hausdorff distance between metric spaces.
\end{itemize}

In terms of heat kernels, the above result has the following counterpart:
\begin{proposition}
Uniformly on compact sets of $\mathbb{R}_{\ge 0} \times \mathbb{R}$.
\[
\lim_{t \to 0} t^2 p_t (\sqrt{t} r, t z) =2\pi^2 h_1(r,z)
\]
\end{proposition}

\begin{proof} Let $K$ be a compact of $\mathbb R_{\ge 0}\times \mathbb{R}$ and $t>0$ sufficiently small so that
 $(\sqrt t r,t z)\in [0, \frac{\pi}{2}]\times [-\pi,\pi]$ forall $(r,z)\in K$.

According to Proposition \ref{representation} we have

$$  t ^2 p_t(\sqrt t r, t z)= \frac{t^{3/2}} {\sqrt{4\pi}} \int_{-\infty}^\infty e^{-\frac{(y+i tz)^2}{4t}} q_t (\cos\sqrt t r \cosh y) dy.
$$ 
The idea is now to use the estimates (\ref{estimate1}) and (\ref{estimate2}) and to study the two integrals:
 
$$J_1(t,r,z)=\int_{\cosh y \leq \frac{1}{\cos\sqrt t r}} e^{- \frac{(y+i tz)^2+ \arccos^2 (\cos\sqrt t r \cosh y)}{4t}}
\frac{\arccos (\cos\sqrt t r \cosh y)}{\sqrt{1-\cos ^2 \sqrt t r \cosh ^2 y  }} dy
$$
and 
$$J_2(t,r,z)=\int_{\cosh y \geq \frac{1}{\cos\sqrt t r}} e^{- \frac{(y+i tz)^2 -\mathrm{arccosh}^2 (\cos\sqrt t r \cosh y)}{4t}}
\frac{\mathrm{arccosh} (\cos\sqrt t r \cosh y)}{\sqrt{\cos ^2 \sqrt t r \cosh ^2 y -1 }} dy.
$$
It is easily seen that for some constant $C>0$, uniformly on $K$,
 $$ |J_1(t,r,z)| \leq C e^{\frac{t z^2}{4}} \sqrt t r.$$

Therefore $J_1(t,r,z)$ goes uniformly to $0$ on $K$.
 
 Let us now turn to the integral $ J_2(t,r,z)$ and let us show that, uniformly,  $ J_2(t,r,z)$ 
converges to $2\pi^2 h_1(r,z)$. 

Let $\e >0$. Let us observe that $|e ^{\frac{iy z}{2}} e^{-\frac{r^2}{4} y \mathrm{cotanh} y}\frac{y}{\sinh y}|$ is less than 
 $y e^{-y}$ for big $y$ and all $r,z$.

 Note also that for all $1< u \leq \cosh(y/2)$,   
 $$e^{-(\frac{y^2- arcosh ^2 u}{4t})} \frac{\mathrm{arcosh} u}{\sqrt{u^2 -1}}\leq e^{-\frac{y^2}{8t}}
 $$
 and for all $\cosh(y/2)\leq u \leq \cosh(y)$
 $$e^{-(\frac{y^2- \mathrm{arcosh} ^2 u}{4t})} \frac{\mathrm{arcosh} u}{\sqrt{u^2 -1}}\leq \frac{y/2}{\sinh{y/2}}.
 $$
 The last three quantities are integrable and do not depend on $r,z$, so we can find  $y_1>0$ so that 
 $$ \int_{ |y|\geq y_1} e ^{\frac{iy z}{2}} e^{-\frac{r^2}{4} y \mathrm{cotanh} y}\frac{y}{\sinh y} dy \leq \e.
 $$
 and
 
 $$\int_{|y| \geq y_1} e^{- \frac{(y+i tz)^2 -\mathrm{arccosh}^2 (\cos\sqrt t r \cosh y)}{4t}}
\frac{\mathrm{arccosh} (\cos\sqrt t r \cosh y)}{\sqrt{\cos ^2 \sqrt t r \cosh ^2 y -1 }} dy \leq \e.
 $$
 Now we study the  behaviour of our integrals for small $y$.
 $|e ^{\frac{iy z}{2}} e^{-\frac{r^2}{4} y \mathrm{cotanh} y}\frac{y}{\sinh y}|$ is less than 
 $1$ for small $y$ and $e^{- \frac{(y+i tz)^2 -arccosh^2 (\cos\sqrt t r \cosh y)}{4t}}
\frac{\mathrm{arccosh} (\cos\sqrt t r \cosh y))}{\sqrt{\cos ^2 \sqrt t r \cosh ^2 y -1 }}$ 
is less than $e^{\frac{t z^2 }{4}}$ for small $y$.
 Thus, as before there exists $0<y_0$ such that 
  
  $$ \int_{ |y|\leq y_0} e ^{\frac{iy z}{2}} e^{-\frac{r^2}{4} y \mathrm{cotanh} y}\frac{y}{\sinh y} dy \leq \e.
 $$
 and
 
 $$\int_{ \mathrm{arccosh}(\frac{1}{\cos \sqrt t r })\leq |y| \leq y_0} e^{- \frac{(y+i tz)^2 -
\mathrm{arccosh}^2 (\cos\sqrt t r \cosh y)}{4t}}
\frac{\mathrm{arccosh} (\cos\sqrt t r \cosh y))}{\sqrt{\cos ^2 \sqrt t r \cosh ^2 y -1 }} dy \leq \e.
 $$
 
 Let $y_0<y<y_1$ and $0<u \leq \cosh y -1$ by the Taylor-Lagrange development formula we have the following equality
 $$\mathrm{arccosh} (\cosh y - u)= y - \frac{1}{\sinh y} u - \frac{ \tilde y}{\sinh ^{3/2} \tilde y} \frac{u^2}{2}.  
 $$ for some $\tilde y \in ]\mathrm{arccosh} (\cosh y - u),  y[$.
 By applying this to $\cos \sqrt t r  \cosh y= \cosh y - t r^2 \cosh y + O (t^2 r^4) \cosh y$, we get
   $$\mathrm{arccosh} (\cos \sqrt t r \cosh y )= y - tr^2 \mathrm{cotanh} y +  O (t^2 r^4) 
(\mathrm{cotanh} y + \cosh ^2 y  \frac{ \tilde y}{\sinh ^{3/2} \tilde y})
 $$ for some $\tilde y \in ]\mathrm{arcosh} (\cos \sqrt  t r \cosh y ),  y[$.
 So 
 $$ \mathrm{arccosh}^2 (\cos \sqrt t r \cosh y )= y^2 - tr^2 y 
\mathrm{cotanh} y +  O (t^2 r^4) ( y \mathrm{cotanh} y +  y \cosh ^2 y  \frac{ \tilde y}{\sinh ^{3/2} \tilde y}).
 $$
 and 
 $$e^{- \frac{(y+i tz)^2 -\mathrm{arccosh}^2 (\cos\sqrt t r \cosh y)}{4t}}= e^{\frac{-iyz}{2}}  e^{- \frac{r^2}{4} y
 \mathrm{cotanh} y}  e^{\frac{tz^2}{4}}
  (1+ O (t^2 r^4) ( y \mathrm{cotanh} y +  y \cosh ^2 y  \frac{ \tilde y}{\sinh ^{3/2} \tilde y}))
 $$
  Finally, using also Taylor Lagrange development formula at order 1 we obtain
  
  $$\frac{\mathrm{arccosh} (\cos\sqrt t r \cosh y)}{\sqrt{\cos ^2 \sqrt t r \cosh ^2 y -1 }}= \frac{y}{\sinh y}
 - \frac{t r^2}{2}\cosh(y) (\frac{1}{\sinh ^2 \hat y} + 2 \frac{\hat y  \cosh \hat y} {\sinh ^3 \hat y} )
  $$ for some $\hat y \in ]\mathrm{arcosh} (\cos \sqrt  t r \cosh y ),  y[$.
  
So finally, we see we can pass  uniformly to the limit under the integral for $y_0\leq |y|\leq y_1$ and obtain our proposition.
  
 
 \end{proof}

This dilation of $\mathbf{SU}(2)$ toward the Heisenberg group can also be seen at the level of differential operators.

Through the map
\begin{eqnarray*}
\mathbf{SU}(2) &\rightarrow &\mathbb{H} \\
\exp(r(\cos \theta X +\sin \theta Y)) \exp {zZ} & \rightarrow & (r,\theta,z)
\end{eqnarray*}
we can see the vector fields $X$, $Y$ and $Z$ of $\mathbf{SU}(2)$ as first order differential operators acting 
on smooth functions on the Heisenberg group with support included in a small enough Carnot Carath\'eodory ball of radius $R$.

Let us now denote by $D$ the dilation vector field on $\mathbb{H}$ given in cylindric coordinates by
\[
D=r \frac{\partial}{\partial r} +2z\frac{\partial}{\partial z}
\] 
For $c\ge 1$ we denote by $X^c$, $Y^c$ and $Z^c$ the dilated vector fields
\[
X^c=\frac{1}{\sqrt{c}} e^{-\frac{1}{2} \ln c D} X e^{ \frac{1}{2} \ln c D} ,
\]
\[
Y^c=\frac{1}{\sqrt{c}} e^{-\frac{1}{2} \ln c D} Y e^{ \frac{1}{2} \ln c D} ,
\]
\[
Z^c=\frac{1}{\sqrt{c}} e^{-\frac{1}{2} \ln c D} Z e^{ \frac{1}{2} \ln c D} .
\]
In the cylindric coordinates of the Heisenberg group, we have

\[
X^c=\cos (-\theta +\frac{2z}{c}) \frac{\partial}{\partial r}+\sin (-\theta + \frac{2z}{c} )
 \left( \sqrt{c}\tan \frac{r}{\sqrt{c}} \frac{\partial}{\partial z}
+\left(\frac{\tan \frac{r}{\sqrt{c}}}{\sqrt{c} }+\frac{1}{\sqrt{c}\tan \frac{r}{\sqrt{c}}}\right) 
 \frac{\partial}{\partial \theta}\right),
\]

\[
Y^c=-\sin(\frac{2z}{c} -\theta) \frac{\partial}{\partial r}+\cos (\frac{2z}{c} -\theta) \left( \sqrt{c} \tan \frac{r}{\sqrt{c}} 
 \frac{\partial}{\partial z}
+\left(\frac{\tan \frac{r}{\sqrt{c}}}{\sqrt{c}} +\frac{1}{\sqrt{c}\tan \frac{r}{\sqrt{c}}}\right)  \frac{\partial}{\partial \theta}\right),
\]
\[
Z^c=\frac{\partial}{\partial z},
\]
so that the dilated vector fields are well-defined on the Carnot-Caratheodory ball with radius $R\sqrt{c}$. 
Consequently, if $f: \mathbb{H} \rightarrow \mathbb{R}$ is a smooth function with compact support, we can speak of $X^cf$, $Y^c f$, and $Z^c f$ as 
soon as the dilation factor  $c$ is big enough. 
For the dilated sublaplacian
\begin{align*}
\mathcal{L}^c& =\frac{1}{c} e^{-\frac{1}{2} \ln c D} \mathcal{L} e^{ \frac{1}{2} \ln c D} \\
  &=(X^c)^2 + (Y^c)^2 \\
  & =\frac{\partial^2}{\partial r^2}+\frac{2}{\sqrt{c}} \mathrm{ cotan } \frac{2r}{\sqrt{c}}\frac{\partial}{\partial r}+\frac{1}{c} \left( 2+\frac{1}{\tan^2 \frac{r}{\sqrt{c}}}
+\tan^2 \frac{r}{\sqrt{c}}  \right) \frac{\partial^2}{\partial \theta^2} +c\tan^2 \frac{r}{\sqrt{c}} \frac{\partial^2}{\partial z^2}+2(1+\tan^2 \frac{2r}{\sqrt{c}})
\frac{\partial^2}{\partial z \partial \theta},
\end{align*}
the same remarks hold true.

With these notations, the \textit{operator} analogue of the convergence of dilated $\mathbf{SU}(2)$ to $\mathbb{H}$ is the following:

\begin{proposition}
If $f:\mathbb{H} \rightarrow \mathbb{R}$ is a smooth function with compact support, then, uniformly, 
$
\lim_{c \to +\infty} X^c f= \tilde{X}f ,
$
$
\lim_{c \to \infty} Y^c f= \tilde{Y} f , 
$
$
\lim_{c \to \infty} Z^c f= \tilde{Z} f  ,
$
$
\lim_{c \to \infty} \mathcal{L}^c f=\tilde{\mathcal{L}} f.
$
\end{proposition}


\section{Gradient bounds for the heat kernel measure}

Let us recall that
\[
\mathcal{L}=X^2+Y^2
\]
with
\[
[X, Y]=2Z, \quad [Y,Z]=2X, \quad [Z,X]=2Y.
\]
In this section, our main goal will be to quantify the regularization property of the semigroup $P_t=e^{t \mathcal{L}}$:
We shall mainly be concerned with bounds for $(XP_tf)^2+(YP_tf)^2$.

\

We shall often make use of the  following notations  (see \cite{livre}, \cite{bakry-st-flour}):
We set for $f,g$ smooth functions,
\[
2\Gamma (f,g)=\mathcal{L}(fg)-f\mathcal{L}g-g\mathcal{L}f
\]
and
\[
2\Gamma_2 (f,g)= \mathcal{L} \Gamma (f,g) -\Gamma (f,\mathcal{L}g)-\Gamma (g,\mathcal{L}f).
\]
In the present setting,
\[
\Gamma (f,f)=(Xf)^2+(Yf)^2
\]
and
\begin{align}\label{Gamma2}
\Gamma_2 (f,f)= (X^2f)^2+(Y^2f)^2+\frac{1}{2} \left( (XY+YX)f \right)^2+2 (Zf)^2 +4 \Gamma (f,f)-4 (Xf)(YZf)+4(Yf)(XZf).
\end{align}

In particular, if $f$ is a smooth function that only depends on the variables $r$ and $z$, we obtain
\begin{align*}
\Gamma (f,f) & = \frac{1}{2}(\mathcal{L}(f^{2})-2 f\mathcal{L}f) \\
 & = \left( \frac{\partial f}{\partial r}\right)^2 +\tan^2 r \left( \frac{\partial f}{\partial z}\right)^2 ,
\end{align*}
and 
\begin{eqnarray*}
\Gamma_2 (f,f) & = &\frac{1}{2}(\mathcal{L}\Gamma(f,f)- 2\Gamma(f,\mathcal{L}f) )\\
  & = & \left( \frac{\partial^2 f}{\partial r^2}\right)^2+
2\tan^2 r \left( \frac{\partial^2 f}{\partial r \partial z}\right)^2
+\tan^4 r \left( \frac{\partial^2 f}{\partial z^2}\right)^2 
  +\frac{2}{\cos^4 r}
\left( \frac{\partial f}{\partial z}\right)^2+\frac{4}{\sin^2 2r} \left( \frac{\partial f}{\partial r}\right)^2\\
 & &
+ \frac{4 \tan r}{\cos^2r}\left( \frac{\partial f}{\partial z}\right)\left( \frac{\partial^2 f}{\partial r \partial z}\right)
- \frac{2 \tan r}{\cos^2r} \left( \frac{\partial f}{\partial r}\right)\left( \frac{\partial^2 f}{\partial z^2}\right) \\
 & =  & \left( \frac{\partial^2 f}{\partial r^2}\right)^2 + \left(  \tan^2 r \frac{\partial^2 f}{\partial z^2} -\frac{2}{\sin 2r}  \frac{\partial f}{\partial r} \right)^2 + 2 \left( \frac{1}{\cos^2 r} \frac{\partial f}{\partial z} +\tan r \frac{\partial^2 f}{\partial r \partial z} \right)^2
\end{eqnarray*}

\subsection{A first gradient bound}

\begin{proposition}
Let $f: \mathbf{SU}(2) \rightarrow \mathbb{R}$ be a smooth function. For $t >0$ and $g \in \mathbf{SU}(2)$,
\[
\Gamma (P_t f , P_t f)(g) \le A(t)
\left( \int_{\mathbf{SU}(2)} f^2 d\mu - \left( \int_{\mathbf{SU}(2)} f d\mu\right)^2 \right)
\]
where
\[
A(t)=-\frac{1}{4} \frac{\partial}{\partial t} \int_{\mathbf{SU}(2)} p^2_t d\mu.
\]

\end{proposition}
\begin{proof}
By left invariance, it is enough to prove this inequality at $g=0$. We can moreover assume that 
$\int_{\mathbf{SU}(2)} f d\mu =0$.
If we denote by $\hat{X}$ and $\hat{Y}$ the right
 invariant vector fields, then we have:
\begin{align*}
\Gamma (P_t f , P_t f)(0)&= (XP_t f)^2(0)+(YP_tf)^2(0) \\
 &=(P_t \hat{X}f)^2(0)+(P_t\hat{Y}f)^2(0) \\
 &=\left( \int_{\mathbf{SU}(2)} p_t (r,z) \hat{X}f (r,z) d\mu \right)^2+
\left( \int_{\mathbf{SU}(2)} p_t (r,z) \hat{Y}f (r,z) d\mu \right)^2 \\
 &=\left( \int_{\mathbf{SU}(2)} \hat{X}p_t (r,z) f (r,z) d\mu \right)^2+
\left( \int_{\mathbf{SU}(2)} \hat{Y}p_t (r,z) f (r,z) d\mu \right)^2
\end{align*}
Now, let us observe that since $p_t$ does not depend on $\theta$, we have
\[
\int_{\mathbf{SU}(2)}  (\hat{X} p_t)^2  d\mu=\int_{\mathbf{SU}(2)}  (\hat{Y} p_t)^2 d\mu
=\frac{1}{2} \int_{\mathbf{SU}(2)}  \Gamma (p_t,p_t) d\mu,
\]
and
\[
\int_{\mathbf{SU}(2)}  \hat{X} p_t \hat{Y} p_t   d\mu=0.
\]
Therefore, from Cauchy-Schwarz inequality, we conclude that:
\[
\Gamma (P_t f , P_t f)(0) \le  \frac{1}{2} \int_{\mathbf{SU}(2)}  \Gamma (p_t,p_t)  d\mu \int_{\mathbf{SU}(2)} f^2 d\mu ,
\]
which is the required inequality because:
\[
\int_{\mathbf{SU}(2)}  \Gamma (p_t,p_t) d\mu=-\int_{\mathbf{SU}(2)}  p_t Lp_t   d\mu
=- \frac{1}{2} \frac{\partial}{\partial t} \int_{\mathbf{SU}(2)} p^2_t  d\mu
\]
\end{proof}
We now study the constant $A(t)$.

\begin{proposition}
We have the following properties:
\begin{itemize}
\item $A$ is decreasing;
\item $A(t) \sim_{t \rightarrow 0} \frac{\pi^2}{32 t^3}$;
\item $A(t) \sim_{t \rightarrow +\infty} 4 e^{-4t}$.
\end{itemize}
\end{proposition}

\begin{proof}
Let us first show that $C$ is decreasing.  We have:
\[
A'(t)=- \int_{\mathbf{SU}(2)} \Gamma_2  ( p_t,  p_t) d\mu.
\]
Since $p_t$ only depends on $(r,z)$, $\Gamma_2 ( p_t, p_t) \ge 0$ and thus $A'(t) \le 0$.

\

We can now observe that, due to the semigroup property,
\[
\int_{\mathbf{SU}(2)} p^2_t d\mu=p_{2t} (0)
\]

But from Proposition \ref{spectral_decomposition} and \ref{asymptotics1},

\begin{align*}
p_t(0)& =\sum_{n=-\infty}^{+\infty}\sum_{k=0}^{+\infty}(2k+\mid n \mid +1)e^{-(4k(k+\mid n \mid +1)+2\mid n \mid  )t} \\
 &=\frac{\pi^2 e^t}{4t^2} \sum_{k \in 
\mathbb{Z}}e^{-\frac{k(k+1)\pi^2}{t}}
\frac{ (2k+1)+2k e^{-\frac{2k\pi^2}{t} }} { \left( 1+e^{-\frac{k\pi^2}{t} } \right)^2},
\end{align*}
which implies the expected result.
\end{proof}

\subsection{Li-Yau type inequality}

We now provide a Li-Yau type estimate for the heat semigroup. The inequality we obtain
 is an improvement in the specific case of $\mathbf{SU}(2)$ of the Cao-Yau gradient estimate 
for subelliptic operators that was obtained in \cite{Cao-Yau}. The idea of the method that is used to prove Theorem \ref{Li-Yau} is due to
 D. Bakry and was given to the authors during personal discussions; It is close to \cite{Bakry-Ledoux}.

We have the following inequality:
\begin{theorem}\label{Li-Yau}
For all $\alpha>2$, for every positive function $f$ and $t >0$,
\[
\Gamma (\ln P_t f) +\frac{t}{\alpha}  (Z \ln P_t f)^2  \leq
\left( \frac{3 \alpha -1}{\alpha -1} -\frac{2 t}{\alpha} \right) \frac{\mathcal{L}P_t f}{P_t f}
+ \frac{t}{\alpha}  - \frac{3\alpha-1}{\alpha-1} + \frac{(3\alpha-1)^2}{\alpha-2}
\frac{1}{t}\]
\end{theorem}

\begin{proof}
We fix a positive function $f$ and $t>0$ and all the following computations are made at a given point $x \in \mathbf{SU}(2)$.

For $0 \le s \le t$, let
\[
\Phi_1 (s)=P_s ((P_{t-s} f) \Gamma (\ln P_{t-s} f))
\]
and
\[
\Phi_2 (s)=P_s ((P_{t-s} f) (Z \ln P_{t-s} f)^2).
\]
Straightforward, but heavy, computations show that
\[
\Phi'_1 (s)=2P_s ((P_{t-s} f) \Gamma_2 (\ln P_{t-s} f))
\]
and
\[
\Phi'_2 (s)=2P_s ((P_{t-s} f) \Gamma(Z \ln P_{t-s} f)).
\]
Now, thanks to the Cauchy-Schwarz inequality, the expression \ref{Gamma2}, shows that for every $\lambda >0$,  and every smooth function $g$,
\[
\Gamma_2 (g) \ge \frac{1}{2} (\mathcal{L}g)^2 +2 (Zg)^2 +\left(4-\frac{2}{\lambda} \right) \Gamma (g)-2\lambda \Gamma (Zg).
\]   
We therefore obtain the following differential inequality
\[
\Phi'_1 (s) \ge P_s ((P_{t-s}f) (\mathcal{L} \ln P_{t-s}f)^2) +4 \Phi_2 (s) +\left(8-\frac{4}{\lambda} \right)\Phi_1 (s) -2 \lambda \Phi_2'(s).
\]
We now have that for every $\gamma \in \mathbb{R}$,
\[
(\mathcal{L} \ln P_{t-s} f)^2 \ge 2 \gamma \mathcal{L} \ln P_{t-s}f - \gamma^2,
\]
and
\[
\mathcal{L} \ln P_{t-s}f =\frac{\mathcal{L} P_{t-s} f}{P_{t-s}f} -\frac{\Gamma (P_{t-s}f)}{(P_{t-s}f)^2 }.
\]
Thus, for every $\lambda >0$ and every $\gamma \in \mathbb{R}$,
\[
\Phi'_1 (s) \ge \left(8-\frac{4}{\lambda} -2 \gamma \right) \Phi_1 (s) +4 \Phi_2 (s) -2 \lambda \Phi_2' (s) +2 \gamma \mathcal{L}P_t f - \gamma^2 P_t f.
\]
Let now $b$ a positive decreasing function on the time interval $[0,t)$. By choosing in the previous inequality
\[
\lambda=-2 \frac{b}{b'}
\]
and
\[
\gamma=\frac{1}{2} \left( 8+2 \frac{b'}{b} + \frac{b''}{b'}\right),
\]
we get
\[
\left( -\frac{1}{4} b' \Phi_1 +b \Phi_2 \right)' \ge -\frac{1}{4} b' \left( \left( 8+2 \frac{b'}{b} + \frac{b''}{b'}\right) \mathcal{L}P_t f -\frac{1}{4}\left( 8+2 \frac{b'}{b} + \frac{b''}{b'}\right)^2 P_t f \right).
\]
Integrating the previous inequality from $0$ to $t$ with the function 
$b(s)=(t-s)^\alpha$, $\alpha >2$,  gives the expected result.
\end{proof}

\begin{remark}
Of course, by the same method, we obtain a Li-Yau type inequality for the heat kernel $p_t (r,z)$ itself.
\end{remark}

\begin{remark}
Interestingly, we can obtain an exponential decay in the previous inequality. Indeed, if we use the function $b(s)=e^{-\frac{8s}{3}} \left( 1-e^{-\frac{8(t-s)}{3\alpha} }\right)^\alpha$, $\alpha >2$, in the 
previous proof, then we obtain that for every $\alpha >2$, $t>0$,
\[
\Gamma (\ln P_t f ) +\frac{3}{2} \left(1-e^{-\frac{8t}{3\alpha} } \right) (Z \ln P_t f )^2
\le
6 \left(-1+\frac{1}{3\alpha} \right)^2 \frac{\alpha}{\alpha-2} \frac{e^{-\frac{16t}{3\alpha}}}{1- e^{-\frac{8t}{3\alpha}}}
-3\left(-1+\frac{1}{3\alpha} \right)\frac{\alpha}{\alpha-1} e^{-\frac{8t}{3\alpha} } \frac{\mathcal{L}P_t f }{P_t f }.
\]
\end{remark}

As a direct corollary of the Li-Yau type inequality of Theorem \ref{Li-Yau}, we classically deduce (by integrating along 
geodesics) the following Harnack type inequality: There exist positive constant $A_1$ and $A_2$ such that for $0<t_1<t_2<1$ and $g_1,g_2 \in \mathbf{SU}(2)$
\begin{align}\label{Harnack}
\frac{p_{t_1} (g_1)}{p_{t_2} (g_2)} \le \left( \frac{t_2}{t_1} \right)^{A_1} 
e^{A_2 \frac{\delta (g_1,g_2)^2}{t_2-t_1}}
\end{align}
where $\delta (g_1,g_2)$ denotes the Carnot-Caratheodory distance from $g_1$ to $g_2$.

As another corollary we can also prove the following global estimate:

\begin{proposition}\label{racine_gamma_log}
There exists a constant $C >0$ such that for  $t \in (0,1)$, $r \in [0,\pi /2]$, $z \in [-\pi,\pi]$,
\[
\sqrt{\Gamma (\ln p_t)(r,z) } \le C \left( \frac{d(r,z)}{t} +\frac{1}{\sqrt{t}} \right),
\]
where $d(r,z)$ denotes the Carnot Carath\'eodory distance from $0$ to the point with cylindric coordinates $(r,\theta, z)$.
\end{proposition}

\begin{proof}
In what follows, we fix $t \in (0,1)$. Let
\[
\phi (s) =P_s ( p_{t-s}  \ln p_{t-s}  )
\]
so that
\[
\phi'(s)= P_s ( p_{t-s}  \Gamma (\ln p_{t-s} ))
\]
and
\[
\phi''(s)=P_s ( p_{t-s} \Gamma_2 (\ln p_{t-s})).
\]
Since $p_t$ only depends on $(r,z)$ we have $\Gamma_2 (\ln p_{t-s}) \ge 0$ and therefore
$\phi''(s) \ge 0$. By integrating the last inequality from 0 to $t/2$, we obtain
\[
\int_0^{\frac{t}{2}} \phi'(s) ds \ge \frac{t}{2} \phi'(0)
\]
that is 
\[
p_t \Gamma (\ln p_t ) \le \frac{2}{t} \left( P_{t/2} \left( p_{t/2}  \ln p_{t/2}  \right) -p_t \ln p_t  \right).
\]
We finally estimate $P_{t/2} \left( p_{t/2}  \ln p_{t/2}  \right) -p_t \ln p_t $ by using first 
$\ln p_{t/2} (r,z) \le \ln p_{t/2} (0)$ and then the  Harnack inequality (\ref{Harnack}):
\[
\ln \frac{ p_{t/2} (0) }{p_t (r,z)} \le C \left( \frac{ d(r,z)^2}{t} +1 \right).
\]
\end{proof}

\subsection{The reverse spectral gap inequality}

As  in the Heisengroup case (see \cite{bakry-baudoin-bonnefont-chafai}), we can easily obtain a reverse Poincare inequality
with a sharp constant for the subelliptic heat kernel measure on $\mathbf{SU}(2)$.

\begin{proposition}
Let $f: \mathbf{SU}(2) \rightarrow \mathbb{R}$ be a smooth function. For $t >0$ and $g \in \mathbf{SU}(2)$,
\[
\Gamma (P_t f , P_t f)(g) \le C(t)
\left( P_t f^2 (g) - (P_t f)^2 (g) \right)
\]
where
\[
C(t)=-\frac{1}{2} \frac{\partial}{\partial t} \int_{\mathbf{SU}(2)} p_t  \ln p_t d\mu.
\]
\end{proposition}

\begin{proof}
By left invariance, it is enough to prove this inequality at $g=0$. If we denote by $\hat{X}$ and $\hat{Y}$ the right
 invariant vector fields, then, as seen before, we have:
\begin{align*}
\Gamma (P_t f , P_t f)(0)=\left( \int_{\mathbf{SU}(2)} \hat{X}p_t (r,z) f (r,z) d\mu \right)^2+
\left( \int_{\mathbf{SU}(2)} \hat{Y}p_t (r,z) f (r,z) d\mu \right)^2
\end{align*}
Since $p_t$ does not depend on $\theta$, we have
\[
\int_{\mathbf{SU}(2)} \frac{ (\hat{X} p_t)^2 }{p_t} d\mu=\int_{\mathbf{SU}(2)} \frac{ (\hat{Y} p_t)^2 }{p_t} d\mu
=\frac{1}{2} \int_{\mathbf{SU}(2)} \frac{ \Gamma (p_t,p_t) }{p_t} d\mu,
\]
and
\[
\int_{\mathbf{SU}(2)} \frac{ \hat{X} p_t \hat{Y} p_t  }{p_t} d\mu=0.
\]
Therefore, from Cauchy-Schwarz inequality, we conclude that:
\[
\Gamma (P_t f , P_t f)(0) \le  \frac{1}{2} \int_{\mathbf{SU}(2)} \frac{ \Gamma (p_t,p_t) }{p_t} d\mu (P_t f^2) (0),
\]
which is the required inequality because:
\[
\int_{\mathbf{SU}(2)} \frac{ \Gamma (p_t,p_t) }{p_t} d\mu=\int_{\mathbf{SU}(2)} \Gamma (\ln p_t,p_t)  d\mu
=-\int_{\mathbf{SU}(2)} \ln p_t Lp_t   d\mu=- \frac{\partial}{\partial t} \int_{\mathbf{SU}(2)} p_t  \ln p_t d\mu
\]
\end{proof}

\begin{remark}
Due to the use of the Cauchy-Schwarz inequality in the previous proof, we see that the previous inequality is sharp.
\end{remark}

We now study the constant
\[
C(t)=-\frac{1}{2} \frac{\partial}{\partial t} \int_{\mathbf{SU}(2)} p_t  \ln p_t d\mu.
\]

\begin{proposition}
We have the following properties:
\begin{itemize}
\item $C$ is decreasing;
\item $C(t)\sim_{t \rightarrow 0} \frac{1}{t}$;
\item $C(t)\sim_{t \rightarrow +\infty} 4 e^{-4t}$.
\end{itemize}
\end{proposition}

\begin{proof}
Let us first show that $C$ is decreasing.
After some computations, we obtain:
\[
C'(t)=- \int_{\mathbf{SU}(2)} \Gamma_2  (\ln p_t, \ln p_t) p_t d\mu.
\]
But now, let us observe that $p_t$ only depends on $(r,z)$. Therefore $\Gamma_2 (\ln p_t,\ln p_t) \ge 0$ and thus $C'(t) \le 0$.

\

We now study $C(t)$ when $t \to 0$. The idea is that, asymptotically when $t \to 0$, the constant $C(t)$ has to behave like 
the best constant of the reverse spectral gap inequality on the Heisenberg group (see the Section 3.4.). From 
\cite{bakry-baudoin-bonnefont-chafai}, this constant is known to be $1/t$.

We have:
\begin{align*}
tC(t)& =\frac{t}{2} \int_{\mathbf{SU}(2)} p_t \Gamma (\ln p_t,\ln p_t)  d\mu \\
 &=\int_{r=0}^{\frac{\pi}{2\sqrt{t}}} \int_{z=-\frac{\pi}{t}}^{\frac{\pi}{t}} t^{5/2} \frac{\sin 2\sqrt{t} r}{4\pi}
p_t(\sqrt{t} r,tz)  \Gamma (\ln p_t,\ln p_t) (\sqrt{t}r,tz) dr dz
\end{align*}
Now, if we denote
\[
\tilde{X}= \cos \theta \frac{\partial }{\partial r} -\frac{\sin
\theta}{r}\frac{\partial }{\partial \theta}- r \sin
\theta \frac{\partial }{\partial z},
\]
\[
\tilde{Y}= \sin \theta \frac{\partial }{\partial r} +\frac{\cos
\theta}{r} \frac{\partial }{\partial \theta} + r \cos
\theta \frac{\partial }{\partial z},
\]
and
\begin{equation*}
h_t(r,z)= \frac{1}{16 \pi^2} \int_{-\infty}^{+\infty} e^{\frac{i \lambda z}{2}} \frac{\lambda}{\sinh \lambda t}
e^{-\frac{r^2}{4}\lambda  \text{cotanh}  \lambda  t }d \lambda,
\end{equation*}
according to the results of Section 3.4., the following convergences hold
\[
\lim_{t \to 0} t^{3/2} \frac{\sin 2\sqrt{t} r}{2\pi}
p_t(\sqrt{t} r,tz)=2\pi h_1 (r,z) r  
\]
\[
\lim_{t \to 0} t \Gamma (\ln p_t , \ln p_t )( \sqrt t r, t z)=(\tilde{X} \ln h_1 )^2 (r,z) +(\tilde{Y} \ln h_1 )^2 (r,z) .
\]
Moreover, thanks to Proposition \ref{racine_gamma_log}, there exists a constant $C>0$ such that
\[
t \Gamma (\ln p_t , \ln p_t )(\sqrt t r,t z) \le C, \quad t \in (0,1).
\]
We can therefore apply a dominated convergence to obtain:
\[
\lim_{t \to 0} t C(t)=\frac{1}{2} \int_{\mathbb{R}^3} h_1 (r,z) \left(
(\tilde{X} \ln h_1 )^2 (r,z) +(\tilde{Y} \ln h_1 )^2 (r,z)\right) r dr d\theta dz.
\]
This last expression is equal to 1, according to \cite{bakry-baudoin-bonnefont-chafai}.

\

We finally turn to the analysis of $C(t)$ when $t \to +\infty$.
For that, we use the expression
\[
C(t)=\frac{1}{2} \int_{\mathbf{SU}(2)} \frac{ \Gamma (p_t,p_t) }{p_t} d\mu
\]
and the spectral decomposition of Proposition \ref{spectral_decomposition} to get that uniformly on $\mathbf{SU}(2)$,
\[
\Gamma (p_t,p_t) \sim_{t \to +\infty} 16 e^{-4t} \Gamma (\cos r \cos z, \cos r \cos z)
\]
Therefore,
\[
C(t) \sim_{t \to \infty} 8 e^{-4t}\int_{\mathbf{SU}(2)} \Gamma (\cos r \cos z, \cos r \cos z) d \mu,
\]
and we compute
\[
\int_{\mathbf{SU}(2)} \Gamma (\cos r \cos z, \cos r \cos z) d \mu=\frac{1}{2},
\]
to conclude.
\end{proof}

\subsection{$L^p$ gradient bounds}

The goal of this section is to prove the following gradient bounds:

\begin{theorem}\label{Driver-Melcher}
Let $p>1$. There exists a constant $C_p>1$ such that for any smooth 
$f:\mathbf{SU}(2) \rightarrow \mathbb{R}$ and any $g\in \mathbf{SU}(2)$
\[
 \sqrt{\Gamma (P_t f, P_tf)  (g)}  \le C_p e^{-2t} \left( P_t \Gamma (f,f)^{\frac{p}{2}} (g)\right)^{\frac{1}{p}}, \quad t \ge 0.
\]
\end{theorem}

\begin{remark}
Let $f(r,\theta,z)=\cos r \cos z$. In that case, $\mathcal{L} f=-2f$ and $\Gamma (f,f)=\sin^2r$. Therefore 
the exponential decay $e^{-2t}$ is optimal and moreover:
\[
\sin r \le C_p P_t (\sin r)^p
\]
which implies, by letting $t \to \infty$,  $C_p \ge \left( 1+\frac{p}{2} \right)^{\frac{1}{p}}$.
\end{remark}

\begin{remark}
We conjecture that the inequality still holds true for $p=1$.
\end{remark}

\subsubsection{Long-time behaviour}

We first study the long-time behaviour $\Gamma (P_t f, P_tf)$. For that we will rely on a commutation between the complex gradient and 
the semigroup $P_t$ (such a type of commutation involving a Folland-Stein type operator has already been used in the Heisenberg group to study
gradient estimates, see \cite{bakry-baudoin-bonnefont-chafai}).

\

The Lie algebra structure relations lead to:
\begin{align}\label{complex_commutation}
(X+iY)L=(L-4iZ+4)(X+iY),
\end{align}
which leads to the \textit{formal} commutation:
\[
(X+iY)P_t =e^{t(L-4iZ+4)}(X+iY).
\]
In what follows we give a precise  analytical sense to the previous commutation.

\begin{remark}
We can observe that the constant that appears in the commutation is positive, which is quite striking because we expect
an exponential decay. Nevertheless, as we will see below $e^{t(L-4iZ)}$ gives a decay $e^{-6t}$ against complex gradients.
\end{remark}

\begin{lemma}
Let $t >0$ and $r \ge 0$. The function
\[
z \rightarrow p_t (r,z)-\frac{1}{\left( 1-\cos r e^{iz-2t}\right)^2}-\frac{1}{ \left(1-\cos r e^{-iz-2t} \right)^2}
\]
admits an analytic continuation on $\left\{ z \in \mathbb{C}, \mid \mathbf{Im} z \mid < -\ln \cos r+6t \right\}$. The function
\[
z \rightarrow p_t (r,z)
\]
is therefore meromorphic on $\left\{ z \in \mathbb{C}, \mid \mathbf{Im} z \mid < -\ln \cos r+6t \right\}$
with double poles at $-i \left( -\ln \cos r +2t \right)$ and $i \left( -\ln \cos r +2t \right)$.
\end{lemma}
\begin{proof}
This is an easy consequence of the spectral decomposition of $p_t$:
\[
p_t (r,z)=\sum_{n=-\infty}^{+\infty}\sum_{k=0}^{+\infty}(2k+\mid n
\mid +1)e^{-(4k(k+\mid n \mid +1)+2\mid n \mid  )t} e^{inz} (\cos
r )^{\mid n \mid}  P_k^{0,\mid n \mid} (\cos 2r ).
\]
\end{proof}

Let us know observe that if $k=0$ and $n \le 0$,
 \[
(X+iY) e^{inz} (\cos r )^{\mid n \mid}  P_k^{0,\mid n \mid}
(\cos 2r )=0.
\]
If, for $t>0$, $r \ge 0$, $z \in \mathbb{C}-\left\{ -i \left( -\ln \cos r +2t \right) \right\}$, 
$\mid \mathbf{Im} z \mid < -\ln \cos r+6t $, we denote
 \[
p_t^* (r,z)=p_t (r,z)-\frac{1}{\left( 1-\cos r e^{-iz-2t}\right)^2},
\]
we have therefore
\[
(X+iY)p_t=(X+iY)p_t^*.
\]
Combining this with (\ref{complex_commutation}) leads to:
\begin{proposition}
If $f: \mathbf{SU}(2) \rightarrow \mathbb{R}$
 is a smooth function, then
\[
 (X +i Y) P_t f (0)  = e^{4t} \int_{\mathbb{H}} p^*_t (r,
z+4it)(X +i Y)f(r, \theta, z) d\mu, \quad t > 0.
\]
\end{proposition}

And, as a corollary:

\begin{corollary}
There exists $t_0>0$ and $A >0$ such that for any smooth 
$f:\mathbf{SU}(2) \rightarrow \mathbb{R}$,
\[
 \sqrt{\Gamma (P_t f, P_tf)  (0)}  \le A e^{-2t}  P_t \sqrt{\Gamma (f,f)} (0),
 \quad t \ge t_0.
\]
\end{corollary}

\begin{proof}
We denote
\[
\Phi (t)=\sup_{r \in [0,\frac{\pi}{2}]} \sup_{z \in [-\pi,\pi]}
\frac{ \mid p^*_t (r,z+4it) \mid }{p_t (r,z) }.
\]
Since, 
\begin{align*}
 p^*_t (r,z+4it)=& \sum_{n=1}^{+\infty} (n+1) e^{-6nt}e^{inz}(\cos r )^{n} \\  & + \sum_{n=-\infty}^{+\infty}\sum_{k=1}^{+\infty}(2k+\mid n
\mid +1)e^{-(4k(k+\mid n \mid +1)+2\mid n \mid  )t} e^{inz}
e^{-4nt} (\cos r )^{\mid n \mid}  P_k^{0,\mid n \mid} (\cos 2r )
\end{align*}
there exists $t_0 >0$ and $A >0$, such that for $t \ge t_0$,
\[
\Phi (t) \le \frac{\mid p^*_t (0) \mid}{1 -\mid 1-p_t(0) \mid }
\le A e^{-6t},
\]
and the proof is complete.
\end{proof}

\begin{remark}
The above function $\Phi (t)$ explodes when $t \rightarrow 0$.
\end{remark}

\begin{remark}
We conjecture that the ratio $\sup_{r,z \notin \Sigma_{R \sqrt{t}}} \frac{p_t (r,z+4it)}{p_t (r,z)}$ is bounded when 
$t \to 0$, where $R$ is big enough and  $\Sigma_{R \sqrt{t}}$ denotes the Carnot Carath\'eodory ball with radius 
$R \sqrt{t}$. By a partition of unity similar to \cite{bakry-baudoin-bonnefont-chafai}, this would imply that Theorem \ref{Driver-Melcher} also holds for $p=1$.
\end{remark}

\subsubsection{Short-time behaviour}

We now conclude the proof of Theorem \ref{Driver-Melcher}, by showing that the
inequality does not explode when $t \rightarrow 0$. We shall use here the
commutation between left-invariant and right invariant vector fields. It relies
on the following lemma:

\begin{lemma}\label{limit}
Let $q >1$. The limit 
\[
\lim_{t \to 0} \int_{\mathbf{SU}(2)} (\sin 2r)^q \Gamma (\ln p_t , \ln
p_t)^{\frac{q}{2}} (r,z) p_t (r,z) d \mu
\]
is finite.
\end{lemma}

\begin{proof}
The proof is similar to the proof of the second point of Proposition 4.9.: By
scaling and a dominated convergence argument based on Proposition 4.6., we
obtain:
 \begin{align*}
 & \lim_{t \to 0} \int_{\mathbf{SU}(2)} (\sin 2r)^q \Gamma (\ln p_t , \ln
p_t)^{\frac{q}{2}} (r,z) p_t (r,z) d \mu \\
 = & 2^q \int_{\mathbb{R}^3}  r ^q h_1 (r,z) \left((\tilde{X} \ln h_1 )^2 (r,z)
+(\tilde{Y} \ln h_1)^2 (r,z) \right)^{q/2}  r dr d\theta dz,
\end{align*}
which is finite, due to known results on the Heisenberg group (see
\cite{bakry-baudoin-bonnefont-chafai}).
\end{proof}

We can now deduce:
\begin{proposition}
Let $p>1$. There exists a constant $A_p>0$ such that for any smooth 
$f:\mathbf{SU}(2) \rightarrow \mathbb{R}$ and any $g\in \mathbf{SU}(2)$
\[
 \sqrt{\Gamma (P_t f, P_tf)  (g)}  \le A_p  \left( P_t \Gamma
(f,f)^{\frac{p}{2}} (g)\right)^{\frac{1}{p}},  t \in (0,1).
\]
\end{proposition}
\begin{proof}
Due to the fact that the right-invariant vector fields $\hat{X}, \hat{Y}$
commute with $\mathcal{L}$, we get
\[
(X P_t f)(0)=(P_t \hat{X} f)(0)
\]
and 
\[
(Y P_t f)(0)=(P_t \hat{Y} f)(0).
\]
Now, $X,Y,Z$ form a basis at each point, there exist therefore smooth functions
such that:
\[
\hat{X}=\Omega_{1,1} X + \Omega_{1,2} Y +\Omega_{1,3} Z
\]
\[
\hat{Y}=\Omega_{2,1} X + \Omega_{2,2} Y +\Omega_{2,3} Z.
\]
By using $[X,Y]=2Z$ and integrating by parts, we obtain
\[
(X P_t f)(0)=\int_{\mathbf{SU} (2)} \left( \Omega_{1,1} p_t +\frac{1}{2} Y (
\Omega_{1,3} p_t) \right) (Xf) +\left( \Omega_{1,2} p_t -\frac{1}{2} X (
\Omega_{1,3} p_t) \right) (Yf) d\mu
\]
and
\[
(Y P_t f)(0)=\int_{\mathbf{SU} (2)} \left( \Omega_{2,1} p_t +\frac{1}{2} Y (
\Omega_{2,3} p_t) \right) (Xf) +\left( \Omega_{2,2} p_t -\frac{1}{2} X (
\Omega_{2,3} p_t) \right) (Yf) d\mu.
\]
We easily compute
\[
\Omega_{1,3}= \sin \theta \sin 2r
\]
and
\[
\Omega_{2,3}=-\cos \theta \sin 2r.
\]
By using H\"older's inequality the expected result follows from  Lemma
\ref{limit}.
\end{proof}

 {\footnotesize %

\end{document}